\documentclass[12pt]{article}
\usepackage{graphicx,psfrag,epsfig, color,float}
\usepackage{amssymb,amsmath,amscd,amsthm}
\usepackage{graphicx,psfrag,epsfig}

\usepackage{graphicx}
\usepackage[active]{srcltx}

\newtheorem{theorem}{Theorem}[section]

\newtheorem{lemma}[theorem]{Lemma}

\def\be{\color{black}}

\date{}

\begin{document}

\date{}
\title{Averaging in the case of multiple invariant measures for the fast system}
\author{
 M. Freidlin\footnote{Dept of Mathematics, University of Maryland,
College Park, MD 20742, mif@math.umd.edu}, L.
Koralov\footnote{Dept of Mathematics, University of Maryland,
College Park, MD 20742, koralov@math.umd.edu}
} \maketitle

\begin{abstract}
 We consider the averaging principle for deterministic or stochastic systems with a fast stochastic component (family of continuous time
Markov chains depending on the state of the system as a parameter). We show that, due to bifurcations in the simplex of invariant probability measures of the chains,
the limiting system should be considered on a graph or on an open book with certain gluing conditions in the vertices of the graph (or on the bifurcation
surface). 
\end{abstract}



\section{Introduction} \label{intro}

Consider the $d$-dimensional continuous stochastic process $z^\varepsilon_t$ satisfying the equation
\begin{equation} \label{e1}
d {z}^\varepsilon_t = v( \xi^\varepsilon_{t}, z^\varepsilon_t) dt + \varkappa d W_t,~~0 < \varepsilon \ll 1.
\end{equation}
We assume that $v$ is sufficiently smooth, and $\xi^\varepsilon_t = \xi_{t/\varepsilon}$, where $\xi_t$ is a stationary process with sufficiently good mixing properties, 
such as a non-degenerate diffusion on 
a compact manifold or a continuous time Markov chain on a finite state space (we consider the latter case in this paper). The Wiener process $W_t$ is independent
of $\xi^\varepsilon_t$. The coefficient $\varkappa$ is non-negative.

 Put $\bar{v}(z) = \mathrm{E} v(\xi_t,z)$. Then (see, for example, \cite{FW} Section 7.2) 
\begin{equation} \label{e2}
z^\varepsilon_t \rightarrow \bar{z}_t,~~~{\rm as}~~\varepsilon \downarrow 0,
\end{equation} \be
(convergence, in distribution, of the processes),
where $\bar{z}_t$ is the solution of the equation
\begin{equation} \label{e3}
d {\bar{z}}_t = \bar{v}(\bar{z}_t) dt + \varkappa d W_t
\end{equation}
with the same initial condition as $z^\varepsilon_t$. 
The convergence of ${z}^\varepsilon_t$ to $\bar{z}_t$ is preserved if the process $\xi_t$ is not stationary but converges with probability one to
a stationary ergodic process $\tilde{\xi}_t$. In this case, $\bar{v}(z) = \mathrm{E} v(z, \tilde{\xi}_t)$.  Moreover, the fast component
$\xi^\varepsilon_{t}$ in (\ref{e1}) can depend on the slow component. In order to illustrate this point, 
let us focus on the case when the fast motion is governed by a continuous time Markov chain
$\Xi^z_t$ on the finite state space $\{1,...,n\}$. The transition rates for the chain $\Xi^z_t$, which depends on the parameter $z \in \mathbb{R}^d$, will be
denoted by  $q_{ij}(z) \geq 0$, $1 \leq i, j \leq n$, $i \neq j$. Intuitively, the slow motion $z^\varepsilon_t$ is governed, at short time scales, by~(\ref{e1})
with $\xi^\varepsilon_t = \xi_{t/\varepsilon}$ replaced by $\Xi^z_{t/\varepsilon}$.
Yet, we cannot simply say that $\Xi^z_{t/\varepsilon}$ is the fast component of the process
since $z$ itself evolves (although slowly) in time. The fast-slow system $X^\varepsilon_t = (\xi^\varepsilon_{t}, z^\varepsilon_{t})$
can be defined constructively (as in Section~\ref{fssy}) or by describing its generator. Namely, for $1 \leq i \leq n$, consider the operators
\[
L_i u(z) = \frac{\varkappa^2}{2} \Delta u(z) + v(i,z) \nabla u(z),
\]
where $u$ is a function defined on $\mathbb{R}^d$. These operators would govern the evolution of the slow component for
the fixed value $i$ of the fast component in the absence of the fast
motion. The second order term, the Laplacian in our case,  could also be a more general operator to allow for more general diffusion in the slow variable. 
To account for the fast component, we define the operator
\[
A^\varepsilon f (i, z) = \frac{1}{\varepsilon} \left(\sum_{j \neq i} q_{ij}(z) (f(j,z) - f(i,z))\right) + L_i f(i,z),
\]
where $f$ is a function on $\{1,...,n\} \times \mathbb{R}^d$. This operator, with the properly specified domain, is the generator of 
the process $X^\varepsilon_t = (\xi^\varepsilon_{t}, z^\varepsilon_{t})$.

 If $q_{ij}(z) > 0$ whenever $i \neq j$, then the process $\Xi^z_t$ has a unique invariant distribution 
$\mu(z) = (\mu_1(z),...,\mu_n(z))$, and (\ref{e2}) holds with $\bar{v}(z) = \sum_{i = 1}^n \mu_i(z) v(i, z)$. Assume now that there is a closed domain $G$
with a smooth boundary such that the chain $\Xi^z_t$ is
ergodic for $z \notin G$ and 
has, say, two ergodic components $R_1 = \{1,...,m\}$ and $R_2 = \{m+1,...,n\}$ for $z \in G$. Thus the transitions between $R_1$
and $R_2$ are impossible while $z^\varepsilon_t \in G$. Then one can expect that, as long as  $z^\varepsilon_t$ remains in $G$, it converges, as $\varepsilon \downarrow 0$, to the solution 
of (\ref{e3}) either with $\bar{v}(z) = \sum_{i \in R_1} \mu_i(z) v(i, z)/\sum_{i \in R_1} \mu_i(z)$ or with 
$\bar{v}(z) = \sum_{i \in R_2} \mu_i(z) v(i, z)/\sum_{i \in R_2} \mu_i(z)$, depending on whether the fast component evolves in $R_1$ or $R_2$. Note that, while
the invariant distribution is not determined uniquely for $z \in G$, the above expressions for $\bar{v}$ are. 

The process $z^\varepsilon_t$ can go from $G$ to $\mathbb{R}^d \setminus G$ and vice versa in finite time. Therefore, in order to define
the limiting process, one should describe the behavior of the process in an infinitesimal neighborhood of $\partial G$.
The novelty of the current work is that, in the presence of multiple invariant measures for the fast process, the limiting motion for
the slow component is (and needs to be) considered on a graph or an open book (if $d > 1$), i.e., a more sophisticated space than the case of non-degenerate 
fast component, where the limiting process lives on the Euclidean space. For simplicity, we'll consider the one-dimensional case, where the
structure of the simplex of invariant probability measures is already non-trivial.

Solutions of the Cauchy problem and of various initial-boundary problems for PDE systems related to the operator $A^\varepsilon$ can be written as expectations of
certain functionals of the process $X^\varepsilon_t = (\xi^\varepsilon_{t}, z^\varepsilon_{t})$. This allows one to calculate the asymptotics of solutions to
those PDE problems using the results for the process ${X}^\varepsilon_t$ and vice versa. One can also apply the probability
results to certain non-linear PDE problems related to the process. For example, certain 
problems for reaction-diffusion systems can be considered in this way (compare
with \cite{F85}, Chapters 5-7). 

Finally, we note that the problem considered in this paper can be viewed as a problem concerning the long-time influence of small perturbations: the process
$\tilde{X}^\varepsilon_t = X^\varepsilon_{\varepsilon t}$ starting at $(i,z)$ can be viewed as  a small perturbation of the process $\tilde{X}_t$ whose first component
is $\Xi^z_t$ starting at $i$ and the second component $z \in \mathbb{R}^d$ does not evolve in time. 

 A general approach to the study of the long-time influence of perturbations (see \cite{Fp1}, \cite{Fp2}) is to consider the projection of $X^\varepsilon_t = \tilde{X}_{t/\varepsilon}$ onto the simplex of invariant 
probability measures of the unperturbed process. In the case when the unperturbed process is $\tilde{X}_t$, the set $M_{\rm erg}$ of the extreme 
points of the simplex (ergodic invariant measures) consists of the measures of the form $\mu(z) \times \delta_z$ (where $z \notin G$ and $\mu(z)$ is the invariant measure
for $\Xi^z_t$) and of the measures of the form $\mu^1(z) \times \delta_z$ and $\mu^2(z) \times \delta_z$ (where $z \in G$ and $\mu^1(z)$, $\mu^2(z)$ are invariant for 
$\Xi^z_t$ on $R_1$ and $R_2$, respectively).  The projection of a point $(i,z)$, $i \in \{1,...,n\}$, $z \in \mathbb{R}^d$, from the phase space of $X^\varepsilon_t$
onto $M_{\rm erg}$ is $\mu(z) \times \delta_z$ (if $z \notin G$), or $\mu^1(z) \times \delta_z$ (if $z \in G$, $i \in R_1$), or $\mu^2(z) \times \delta_z$
(if $z \in G$, $i \in R_2$).  Note that $M_{\rm erg}$
can be parametrized by the set of pairs $(l,z)$, $z \in \mathbb{R}^d$, $l \in \{1, 2\}$ if $z \in G$ and $l = 0$ if $z \notin G$.
The main result of the paper is that the projection of $X^\varepsilon_t$ onto $M_{\rm erg}$ converges to a Markov process on $M_{\rm erg}$.  \be


\section{The fast-slow system} \label{fssy}

In this section, we'll introduce the fast-slow system $X^\varepsilon_t = (\xi^\varepsilon_{t}, z^\varepsilon_{t})$. (Sometimes we'll write $X^{x,\varepsilon}_t$ to
indicate the dependence on the initial position $x$). The fast component $\xi^\varepsilon_{t}$ evolves as a Markov chain, 
whose transition rates depend on the slow variable. The slow component $z^\varepsilon_{t}$ solves an ODE or an SDE with the right hand side that depends on the
fast variable.
 Namely,  let $q_{ij}(z) \geq 0$, $1 \leq i, j \leq n$, $i \neq j$ be a family of transition rates
for a Markov chain  $\Xi^z_t$ that depends on the parameter $z \in \mathbb{R}$. Each of the functions $q_{ij}(z)$ is assumed to be continuous. 

We assume that  $\Xi^z_t$ is ergodic for each $z < 0$, while there are two ergodic classes $R_1 = \{1,...,m\}$ and $R_2 = \{m+1,...,n\}$ for $z > 0$. 
More precisely, we us assume that $q_{ij}(z) > 0$  
for $z < 0$, while, for $z \geq 0$, $q_{ij}(z) > 0$ if and only if $i, j \in R_1$ or $i, j \in R_2$. 

Moreover, we assume that $q_{ij}(z) - q_{ij}(0)$, $i \neq j$, degenerate at the same rate as $z \uparrow 0$, namely, there are
positive constants $\overline{q}_{ij}$, a function $\varphi: (-\infty, 0) \rightarrow (0, \infty)$ with $\lim_{z \uparrow 0} \varphi(z) = 0$ and
functions $\beta_{ij}: (-\infty, 0) \rightarrow \mathbb{R}$ with $\lim_{z \uparrow 0} \beta_{ij}(z) = 0$ such that
\[
q_{ij}(z) - q_{ij}(0) = \overline{q}_{ij} \varphi(z) (1 + \beta_{ij}(z)),~~~~z < 0,~~~i \neq j.
\]
Let $\mu_i(z)$, $1 \leq i \leq n$, $z \in  \mathbb{R}$,
be the invariant distribution of the Markov chain $\Xi^z_t$. This is not determined uniquely for $z \geq 0$ since there are two ergodic classes for the
Markov chain. However, under the above assumptions, there are limits
$\pi_i = \lim_{z \uparrow 0} \mu_i(z)$, and, for $z \geq 0$, we select the unique invariant distribution such that $\mu_i(z)$ are continuous functions on $\mathbb{R}$.
Define 
\[
\overline{\pi}_1 = \sum_{i \in R_1} \pi_i,~~~\overline{\pi}_2 = \sum_{i \in R_2} \pi_i.
\] 
Let $v(i,z)$, $1 \leq i \leq n$, $z \in \mathbb{R}$, be Lipschitz-continuous in $z$ for each $i$. Define
\[
\overline{v}_0(z) = \sum_{i=1}^n v(i, z) \mu_i(z),~~z < 0,
\]
\[\overline{v}_1(z) = \frac{1}{\overline{\pi}_1} \sum_{i \in R_1} v(i, z) \mu_i(z),~~
\overline{v}_2(z) =  \frac{1}{\overline{\pi}_2} \sum_{i \in R_2} v(i, z) \mu_i(z),~~z \geq 0.
\]
 We'll assume that $v(i,z) > 0$ for each $(i,z)$ (this assumption is not required if there is
diffusion in the slow variable (case $\varkappa = 1$ below)). Let us make a simplifying assumption about the behavior of the coefficients at infinity. Namely,  we will assume that there is $C > 0$ such that
$q_{ij}(z) = q_{ij}^l$ for $z \leq -C$ and $q_{ij}(z) = q_{ij}^r$ for $z \geq C$, where $q_{ij}^l, q_{ij}^r$ do not depend on $z$. 
Moreover, let us assume that $v(i,z) = \overline{v}^\infty$ for some  $\overline{v}^\infty$  for all $1 \leq i \leq n$, $|z| \geq C$. These assumptions can be relaxed significantly, 
however, this will not concern us  since we would like to focus on the behavior of the process near $z = 0$. The slow component $z^\varepsilon_{t}$ is assumed to be continuous and to satisfy
\[
{d z^\varepsilon_{t}} = v(\xi^\varepsilon_{t},z^\varepsilon_{t})dt + \varkappa d W_t
\]
at the points of continuity of $\xi^\varepsilon_{t}$. Here $\varkappa = 0$ or $\varkappa = 1$ (we'll consider two cases resulting in two different types of the limiting behavior). The fast component, intuitively, evolves as the Markov chain $\Xi^z_t$ (with $z = z^\varepsilon_t$), sped up by the factor $1/\varepsilon$. However, since $z$ itself evolves in time, we need a more formal
definition of the process $X^{x,\varepsilon}_t = (\xi^{x,\varepsilon}_{t}, z^{x,\varepsilon}_{t})$. Namely, the process 
starts at $x  = (i,z) \in M$, and moves along the
$z$-axis during a random time interval $[0,\sigma)$. For $t \in [0,\sigma)$, $z^{x,\varepsilon}_{t}$ solves
\[
{d z_{t}} =  v(i,z_{t}) dt + \varkappa d W_t.
\]
At a random time $\sigma$, the process $X^{x,\varepsilon}_t$ jumps to a random location $(j, z_{\sigma})$. The distribution of $\sigma$ is determined 
as follows. Let $Q_i(z) = \sum_{j \neq i} q_{ij}(z)$ and $r(t) = \varepsilon^{-1} \int_0^t Q_i(z_s) ds$. Then $\sigma \geq 0$ is such that the distribution of $r(\sigma)$ is exponential with parameter one. Given a value of $\sigma$, the probability that $X^{x,\varepsilon}_t$ jumps to $(j, z_{\sigma})$ is $q_{ij}(\sigma) / Q_i(\sigma)$. 

Having identified the location of the process at time $\sigma$, we treat it as a new starting point, and select a new (random) time interval for the 
jump-free motion of the process independently of the past. The construction then continues inductively. It is clear that the process just described
is the RCLL Markov process. 

The process $X^{x,\varepsilon}_t = (\xi^{x,\varepsilon}_{t}, z^{x,\varepsilon}_{t})$ could be defined, equivalently,
through its generator, using the Hille-Yosida theorem. We discuss the Hille-Yosida theorem and the generator of $X^{x,\varepsilon}_t$ next, since, in any case, a 
similar construction will be used to define the limiting process when $\varkappa = 1$.

Let $M$ be a separable locally compact metric space,  $C_0(M)$ be the space of continuous functions on $M$ that tend to zero at infinity (can be made arbitrarily 
close to zero outside a sufficiently large compact). The space $C_0(M)$ is endowed with the supremum norm. Let $P(t,x,B)$ be a Markov transition function (a priori not assumed to be conservative) on $M$.
For $f \in C_0(M)$, let 
\[
(T_t f)(x) = \int_{\mathbb{T}^d} f(x')P(t,x,dx'),~~~t \geq 0.
\]
We'll say that $P$ satisfies condition $C_0$ if $T_t f \in C_0(M)$ for each $f \in C_0(M)$. Recall that $P$ is said to
be stochastically continuous if $\lim_{t \downarrow 0} P(t,x,U) = 1$ for each open neighborhood $U$ of $x$.

\begin{theorem} {\rm [Hille-Yosida] (\cite{KPST}, page 365).}
 Suppose that a linear operator $A$ on $C_0(M)$ has the following properties:

(a) The domain $\mathcal{D}(A)$ is dense in $C_0(M)$;

(b) If $f \in \mathcal{D}(A)$, $f(x_0) \geq 0$ and $f(x_0) \geq f(x)$ for all $x \in M$, then $A f(x_0) \leq 0$.  

(c) For every $\psi \in C_0(M)$, and every $\lambda > 0$, there exists a solution $f \in \mathcal{D}(A)$ of the equation $\lambda f - A f = \psi$.

Then the operator $A$ is the infinitesimal generator of a semi-group $T_t$, $t \geq 0$, on $C_0(M)$ that 
is defined by a stochastically continuous Markov transition function satisfying condition $C_0$. The transition function with 
such properties is determined uniquely.
\end{theorem}

The Hille-Yosida theorem can be applied to the space $M = \{1,...,n\} \times \mathbb{R}$. Let us define the linear operator $A^\varepsilon$ in $C_0(M)$. 
In the case $\varkappa = 0$,  the  domain of $A^\varepsilon$, denoted by
$\mathcal{D}(A^\varepsilon)$, consists of all functions $f \in C_0(M)$ such that $f'(i,\cdot) \in C_0(\mathbb{R})$ for each $i$. 
For $f \in \mathcal{D}(A^\varepsilon)$, we define 
\[
A^\varepsilon f (i, z) = \frac{1}{\varepsilon} \left(\sum_{j \neq i} q_{ij}(z) f(j,z) - Q_i( z) f(i,z)\right) + v(i,z) f'(i, z).
\]


In the case when $\varkappa = 1$, the  domain of $A^\varepsilon$ 
consists of all functions $f \in C_0(M)$ such that $\frac{1}{2} f''(i,\cdot) + v(i,\cdot) f'(i,\cdot) \in C_0(\mathbb{R})$ for each $i$. 
For $f \in \mathcal{D}(A^\varepsilon)$, we define 
\[
A^\varepsilon f (i, z) = \frac{1}{\varepsilon} \left(\sum_{j \neq i} q_{ij}(z) f(j,z) - Q_i( z) f(i,z)\right) + \frac{1}{2} f''(i,z) + v(i,z) f'(i, z).
\]

In both cases, it is possible to show that the conditions of the Hille Yosida theorem are satisfied. (We skip details since, in any case, the process was
already defined constructively.) 
%
%
Let $P^\varepsilon(t,x,dx')$ be the corresponding Markov transition function, and  $T^\varepsilon_t$, $t \geq 0$, be the corresponding semi-group on $C_0(M)$.
Take a sequence of functions $f_n \in  \mathcal{D}(A^\varepsilon)$ with values in $[0,1]$  with compact support such that $f_n(i,z) = 1$ for $|z| \leq n$,  
$\|A^\varepsilon f_n\|_{C_0} \leq 1/n$.
 The existence of
such a sequence is easily justified once we recall that the coefficients of $A^\varepsilon$ are constant for sufficiently large $|z|$.   

Since ${A^\varepsilon}$ is the infinitesimal generator of the semi-group $T^\varepsilon_t$, we 
have (see Theorem I.1 of \cite{Mandl}), for $f \in \mathcal{D}(A^\varepsilon)$, 
\begin{equation} \label{sequ}
T^\varepsilon_t f - f = \int_0^t T^\varepsilon_s A^\varepsilon f ds.
\end{equation}
Therefore,
\[
T^\varepsilon_t f_n (x) - f_n(x) = \int_0^t T^\varepsilon_s A^\varepsilon f_n (x) ds \rightarrow 0~~~{\rm as}~~n \rightarrow \infty, 
\]
which implies that $T^\varepsilon_t f_n (x) \rightarrow 1$, and therefore $P^\varepsilon(t,x,\cdot)$ is a probability measure. 
Let $X^{x, \varepsilon}_t = (\xi^{x,\varepsilon}_{t},z^{x,\varepsilon}_{t})$, $x = (i,z)\in M$, be the corresponding Markov family. 
A modification of $X^{x, \varepsilon}_t$ can be chosen with
trajectories that are right continuous and have left limits (\cite{KPST}, page 348).  

Rewrite (\ref{sequ}) as
\[
\mathrm{E} f (X^{x,\varepsilon}_t) - f(x) = \mathrm{E} \int_0^t (A^\varepsilon f) (X^{x,\varepsilon}_s) ds. 
\]
Since $X^{x,\varepsilon}_t$ is a RCLL Markov process with continuous trajectories, for each $x \in M$, the process 
$f (X^{x,\varepsilon}_t) - f(x) - \int_0^t (A^\varepsilon f) (X^{x,\varepsilon}_s) ds$ is a RCLL martingale, and, for each stopping time $\tau$ 
with $\mathrm{E} \tau < \infty$, we get 
\begin{equation} \label{sttime}
\mathrm{E} f (X^{x,\varepsilon}_\tau) - f(x) = \mathrm{E} \int_0^\tau (A^\varepsilon f) (X^{x,\varepsilon}_s) ds.
\end{equation}

Recall that we earlier defined the process $X^{x,\varepsilon}_t$ constructively, without referring to the Hille-Yosida theorem. It is
easily verified directly that the generator of this process coincides with $A^\varepsilon$ on $\mathcal{D}(A^\varepsilon)$. 
The  Markov transition function of the process is stochastically continuous and satisfies condition $C_0$. At the same time, by (\ref{sequ}),
the semigroup is defined uniquely by the values of the generator on a dense set, and thus the generator of the constructively defined process is
$A^\varepsilon$ (rather than a non-trivial extension).

\section{The limiting process} \label{limpro}

Let us describe the appropriate space and the limiting process on it for the fast-slow system $X^{x,\varepsilon}_t = (\xi^{x,\varepsilon}_{t}, z^{x,\varepsilon}_{t})$.
Let $I_0 = (-\infty, 0]$, $I_1 = \{1\} \times [0,\infty)$, $I_2 = \{2\} \times [0, \infty)$. These are three half-lines, with $I_1$ and $I_2$ distinguished by
a label. We'll identify the ends of $I_0$, $I_1$, and $I_2$, thus obtaining a graph, denoted by $S$, with three semi-infinite edges with the common vertex, which will be denoted $O$. 
 Each point $y = (l, z) \in S$ is determined by the label of the edge $l \in \{0,1,2\}$ and the coordinate $z$, 
where $z \in (-\infty, 0]$ for $l = 0$ and $z \in [0, \infty)$ for $l =1,2$.  

First, consider the case when there is no diffusion in the slow variable ($\varkappa = 0$). The process $Y^y_t$ starting at $y = (l,z) \in S$ will move deterministically with the variable speed $\overline{v}_0$ on $I_0$, 
$\overline{v}_1$ on $I_1$, and $\overline{v}_2$ on $I_2$. For $y \in I_0$, we still need to describe the behavior of $Y^y_t$ once the process reaches $O$. 
The behavior at $O$ is random, the process proceeds to $I_1$ and $I_2$ with probabilities
\[
p_1 =  \frac{ \sum_{i \in R_1} \pi_i v_i  }{
 \sum_{i \in R_1} \pi_i v_i  +  \sum_{i \in R_2} \pi_i v_i  }~~~~~{\rm and}~~~~~p_2 =  \frac{ \sum_{i \in R_2} \pi_i v_i  }{
 \sum_{i \in R_1} \pi_i v_i  +  \sum_{i \in R_2} \pi_i v_i  }
\]
respectively, where $v_i = v(i,0)$.

Next, consider the case with diffusion ($\varkappa = 1$). The process $Y^y_t$ is a diffusion inside each of the edges. However, a gluing condition is
needed to describe the behavior of the process once it reaches the vertex. Thus, it is most convenient to define the process via its generator.
The  domain of $A$, denoted by
$\mathcal{D}(A)$, consists of all functions $f \in C_0(S)$ such that: 

(a) $\frac{1}{2} f''(l,\cdot) + \overline{v}_l(\cdot) f'(l,\cdot) \in C_0(S)$, i.e., the differential operator can be applied to $f$ inside
each of the edges, and the resulting function can be extended to the vertex $O$, so that it becomes an element of $C_0(S)$.

(b) There are one-sided derivatives $f'(l,0)$ and 
\begin{equation} \label{relder}
f'(0,0) = \overline{\pi}_1 f'(1,0) + \overline{\pi}_2 f'(2,0).
\end{equation} 
\\
It is not difficult to verify that the conditions of the Hille-Yosida theorem are satisfied and that the resulting Markov transition function,
denoted by $P(t,x,B)$, is a probability measure, as a function of $B$. Let $Y^{y}_t$, $y \in S$, be the corresponding Markov family and $T_t$ be the
corresponding semigroup. 
In order to show that a modification with continuous trajectories exists, it is enough to check that $\lim_{t \downarrow 0} P(t,x,B)/t = 0$ for each closed set $B$ that doesn't contain $x$ (Theorem I.5 of \cite{Mandl}, see also \cite{Dyn}). Let  $f \in \mathcal{D}(A)$  be a non-negative function that is equal to one on $B$ and whose support doesn't contain $x$. Then
\[
\lim_{t \downarrow 0} \frac{P(t,x,B)}{t} \leq \lim_{t \downarrow 0} \frac{(T_t f)(x) - f(x)}{t} = A f(x) = 0,
\]
as required. Thus $Y^{y}_t$ can be assumed to have continuous trajectories. 

\be

\section{A lemma on convergence of processes}

The next lemma can be used  to show convergence of families of parameter-dependent processes.
We formulate it in a general setting. Consider a metric space $M$ and a 
Markov family $X^{x, \varepsilon}_t$, $x \in M$, of processes that depend on a 
parameter $\varepsilon > 0$. We also consider a continuous mapping $h: M \rightarrow S$  
from  $M$ to a  locally compact  separable  metric space 
$S$ and define the processes $Y^{x, \varepsilon}_t = h(X^{x, \varepsilon}_t)$, 
$x \in M$, $\varepsilon > 0$. 

The motivation to introduce the latter family of processes comes from our desire to study the limiting behavior of  $X^{x, \varepsilon}_t$, as $\varepsilon \downarrow 0$.
However, the space $M$ is too large for our purposes, i.e., the natural state space for the limiting process consists of equivalence classes
in $M$ rather than of individual points. 
Thus, $Y^{x, \varepsilon}_t$ will capture reduced dynamics, where meaningful limiting behavior can be observed.  

Note that while convergence to Markov processes on $S$ as $\varepsilon \downarrow 0$ will be established, the processes $Y^{x, \varepsilon}_t$
need not be Markov for fixed $\varepsilon > 0$. The main point of the lemma is that,  in order to demonstrate the convergence of 
$Y^{x, \varepsilon}_t$ to a limiting process, it is sufficient to check that for small $\varepsilon$ the processes
nearly satisfy the relation (\ref{mprob}), which is similar to the martingale problem but with the ordinary exptectation rather than the conditional 
expectation.

\begin{lemma} \label{cnvl} Let $h: M \rightarrow S$ be a continuous mapping from a metric space $M$ to a locally compact separable metric space $S$. Let $X^{x, \varepsilon}_t$, $x \in M$, be a  Markov family on $M$ that depends on a parameter $\varepsilon > 0$. Suppose that the processes $Y^{x, \varepsilon}_t = h(X^{x, \varepsilon}_t)$, 
$x \in M$, $\varepsilon > 0$, have
continuous trajectories. Let $Y^y_t$, $y \in S$, be a Markov family on $S$ with continuous trajectories whose semigroup 
$T_t$, $t \geq 0$,  preserves the space $C_0(S)$. 
(This, together with the continuity of trajectories, implies that $T_t$ is a Feller semi-group, i.e.,  $T_t f$, viewed as a function of $t$, is a right-continuous
from $[0,\infty)$ to $C_0(S)$ for each $f$.) Let $A: \mathcal{D}(A) \rightarrow C_0(S)$ denote the infinitesimal generator of this family,  where $\mathcal{D}(A)$ is the domain of the generator.
Let $\Psi$ be dense linear subspace  of $C_0(S)$ and $\mathcal{D}$ be a linear
subspace of $\mathcal{D}(A)$, and
suppose that $\Psi$ and $\mathcal{D}$ have the following properties:

(1)   There is   $\lambda > 0$  such that  for  each $f \in \Psi$ the equation $\lambda F -  A F = f$ has a solution~$F \in \mathcal{D}$. 

(2) For each $T > 0$, each $f \in \mathcal{D}$, and each compact $K \subseteq S$,
  \begin{equation} \label{mprob}
	\lim_{\varepsilon \downarrow 0} \mathrm{E} (f( Y^{x, \varepsilon}_T ) - f (  Y^{x, \varepsilon}_0) - \int_0^T A f (  Y^{x, \varepsilon}_t) dt) =  0,
	\end{equation}
uniformly in $x \in h^{-1} (K)$. 
Suppose that the family of measures on $C([0, \infty), S)$ induced by the processes $Y^{x, \varepsilon}_t$, $\varepsilon > 0$, is tight for each $x \in M$.
 
Then, for each $x \in M$, the measures induced by the processes $Y^{x, \varepsilon}_t$ converge weakly, as $\varepsilon \downarrow 0$, to
the measure induced by the process $Y^{h(x)}_t$. 
\end{lemma}
\proof Fix $x \in M$. Since the family of measures on $C([0, \infty), S)$ induced by the processes $Y^{x, \varepsilon}_t$, $\varepsilon > 0$, is tight, we can find a process $Z^x_t$ with continuous trajectories and a sequence
$\varepsilon_n \downarrow 0$ such that  $Y^{x, \varepsilon_n}_t$ converge to $Z^x_t$ in distribution as $n \rightarrow \infty$. The desired result will immediately follow if we demonstrate that the distribution of $Z^x_t$ 
coincides with the distribution of $Y^{h(x)}_t$ (and thus does not depend
on the choice of the sequence $\varepsilon_n$). We will show that $Z^x_t$ is a solution of the martingale problem for $(A|_\mathcal{D}, h(x))$, i.e.,  for each $T_2 > T_1 \geq 0$ and $f \in \mathcal{D}$,
  \begin{equation} \label{rmpro}
\mathrm{E} (f( Z^{x}_{T_2} ) - f (  Z^{x}_{T_1}) - \int_{T_1}^{T_2} A f (  Z^{x}_t) dt| {\mathcal{F}}_{T_1}^{Z^x}) =  0,~~~~~Z^x_0 = h(x). 
	\end{equation}
First, however, let us discuss the uniqueness for solutions of the martingale problem. 	
We claim that:

(a) $\mathcal{D}$ is dense in $C_0(S)$.

(b) ${\rm Range}(\lambda - A|_\mathcal{D})$ is dense in $C_0(S)$.

(c) For each pair of measures $\mu_1$, $\mu_2$ on $S$, the equality $\int_S f d\mu_1 = \int_S f d \mu_2$ for all $f \in C_0(S)$ implies that
$\mu_1 = \mu_2$.

 To demonstrate (a), 
take an arbitrary $\delta > 0$ and $F_0 \in \mathcal{D}(A)$. Let $g_0 = \lambda F_0 - A F_0$, and take $g' \in \Psi$ such that $\|g' - g_0\| \leq \lambda \delta$. Let $F' \in \mathcal{D}$
be such that $\lambda F' - A F' = g'$. Then, since $A$ is the generator of a strongly continuous semigroup on $C_0(S)$, from the Hille-Yosida theorem  
it follows that $\|F' - F_0\| \leq \|g' - g_0\|/\lambda \leq \delta$. This implies (a) since  $\mathcal{D}(A)$ is dense in $C_0(S)$. Note that (b) follows from the existence of a solution $F \in \mathcal{D}$ to $\lambda F -  A F = f \in \Psi$ and the density of $\Psi$, while (c) is obvious.  The validity of (a)-(c) is enough to conclude that the distribution on $C([0, \infty), S)$ of a process with continuous paths satisfying (\ref{rmpro}) is uniquely determined (Theorem 4.1, Chapter 4 in \cite{EK86}).

Note that (\ref{rmpro}) is satisfied if $Z^x_t$ is replaced by $Y^{h(x)}_t$ since $\mathcal{D} \subseteq \mathcal{D}(A)$ and $A$ the the generator of the family $Y^y_t$, $y \in S$. 
Therefore, $Z^x_t$ and $Y^{h(x)}_t$ have the same distribution if (\ref{rmpro}) holds. It remains to prove (\ref{rmpro}).

Note that $Z^x_t$ is a solution of the martingale problem for $(A|_\mathcal{D}, h(x))$ if and only if 
\[
\mathrm{E} \left( (\prod_{i=1}^k  g_i(Z^x_{t_i})) (f( Z^{x}_{T_2} ) - f (  Z^{x}_{T_1}) - \int_{T_1}^{T_2} A f (  Z^{x}_t) dt) \right) =  0,~~~~~Z^x_0 = h(x),
\]
whenever $f \in \mathcal{D}$, $0 \leq t_1 < ... < t_k \leq T_1$, and $g_1,...,g_k \in C_0(S)$. Since $Y^{x, \varepsilon_n}_t = h(X^{x, \varepsilon_n}_t)$ converge to $Z^x_t$ in distribution, we have
\[
\mathrm{E} \left( (\prod_{i=1}^k  g_i(Z^x_{t_i})) (f( Z^{x}_{T_2} ) - f (  Z^{x}_{T_1}) - \int_{T_1}^{T_2} A f (  Z^{x}_t) dt) \right) =
\]
\[
\lim_{n \rightarrow \infty} \mathrm{E} \left( (\prod_{i=1}^k  g_i(h(X^{x,\varepsilon_n}_{t_i}))) (f( h(X^{x,\varepsilon_n}_{T_2}) ) - f (  h(X^{x,\varepsilon_n}_{T_1})) - \int_{T_1}^{T_2} A f (  h(X^{x,\varepsilon_n}_t)) dt) \right) = 
\]
\[
\lim_{n \rightarrow \infty} \mathrm{E} \left( (\prod_{i=1}^k  g_i(h(X^{x,\varepsilon_n}_{t_i}))) \mathrm{E}(f( h(X^{x,\varepsilon_n}_{T_2}) ) - f (  h(X^{x,\varepsilon_n}_{T_1})) - \int_{T_1}^{T_2} A f (  h(X^{x,\varepsilon_n}_t)) dt| \mathcal{F}_{T_1}^{X^{x,\varepsilon_n}}) \right).
\]
By the Markov property of the  family  $X^{x, \varepsilon_n}_t$, 
\[
\mathrm{E}(f( h(X^{x,\varepsilon_n}_{T_2}) ) - f (  h(X^{x,\varepsilon_n}_{T_1})) - \int_{T_1}^{T_2} A f (  h(X^{x,\varepsilon_n}_t)) dt| \mathcal{F}_{T_1}^{X^{x,\varepsilon_n}}) = 
\]
\[
\mathrm{E}(f( h(X^{x',\varepsilon_n}_{T_2-T_1}) ) - f (  h(X^{x',\varepsilon_n}_{0})) - \int_{0}^{T_2-T_1} A f (  h(X^{x',\varepsilon_n}_t)) dt) |_{x' = X^{x,\varepsilon_n}_{T_1}}~, 
\]
which tends to zero in distribution, as follows from (\ref{mprob}) and   from the tightness of the sequence of random variables $X^{x,\varepsilon_n}_{T_1}$.
Therefore, using the boundedness of $f$, $Af$, and $g_1,...,g_k$, we conclude that 
\[
\mathrm{E} \left( (\prod_{i=1}^k  g_i(Z^x_{t_i})) (f( Z^{x}_{T_2} ) - f (  Z^{x}_{T_1}) - \int_{T_1}^{T_2} A f (  Z^{x}_t) dt) \right) = 0.
\]
Finally, $Z^x_0 = h(x)$ since $Y^{x, \varepsilon_n}_0 = h(X^{x,\varepsilon_n}_0) = h(x)$ for all $n$. \qed

\section{Convergence of the fast-slow process}
\subsection{The case with no diffusion}

Consider first a simplified version of the problem: assume that the fast-slow system $X^{x,\varepsilon}_t = (\xi^{x,\varepsilon}_{t}, z^{x,\varepsilon}_{t})$
is defined as in Section~\ref{fssy}, but  $q_{ij}(z) > 0$ for $i \neq j$ 
(and thus $\Xi^z_t$ is ergodic) for each $z \in \mathbb{R}$.
In this case, the fast Markov chain has a unique invariant distribution, which will be denoted by $\mu_i(z) $, $1 \leq i \leq n$, for each $z \in  \mathbb{R}$.
 Define $Y^y_t$, $y \in \mathbb{R}$, to be the deterministic motion on the real line
with the velocity $\overline{v}(y) = \sum_{i=1}^n v(i, z) \mu_i(z)$, $z \in \mathbb{R}$. The domain $\mathcal{D}(A)$ of its generator $A$ consists of all functions $f \in C_0(\mathbb{R})$ such that $f' 
\in C_0(\mathbb{R})$, while $A f (y)  = \overline{v}(y) f'(y)$. Let $h: M \rightarrow \mathbb{R}$ be the projection $h(i,z) = z$. 
The following theorem is a standard averaging result. 
\begin{theorem} \label{saba} Suppose that $q_{ij} > 0$ for $i \neq j$, $z \in \mathbb{R}$.
For each $x \in M$, the measures induced by the processes $Y^{x, \varepsilon}_t = h(X^{x, \varepsilon}_t)$ on 
$\mathbb{R}$ converge weakly, as $\varepsilon \downarrow 0$, to
the measure induced by the process $Y^{h(x)}_t$. 
\end{theorem}
\proof
We apply Lemma~\ref{cnvl} with $S = \mathbb{R}$, $\Psi=\mathcal{D} = \mathcal{D}(A)$. Thus we need to justify (\ref{mprob}) for $f \in \mathcal{D}(A)$.
Define $\tilde{f}(i,z) = f(z)$, $1 \leq i \leq n$. Using (\ref{sttime}) (which is still valid in this simplified case) applied to $\tilde{f}$ with $\tau = T$, we can write
\[
 \mathrm{E} \left(f( Y^{x, \varepsilon}_T ) - f (  Y^{x, \varepsilon}_0) - \int_0^T A f (  Y^{x, \varepsilon}_t) dt\right) =
\]
\[
\mathrm{E} \left(f( Y^{x, \varepsilon}_T ) - f (  Y^{x, \varepsilon}_0) - \int_0^T A f (  Y^{x, \varepsilon}_t) dt\right)  - 
\mathrm{E} \left(\tilde{f}( X^{x, \varepsilon}_T ) - \tilde{f} (x) - \int_0^T A^\varepsilon \tilde{f} (  X^{x, \varepsilon}_t) dt\right) =
\]
\[
\mathrm{E}  \int_0^T \left( A^\varepsilon \tilde{f} (  X^{x, \varepsilon}_t) -  A f (  Y^{x, \varepsilon}_t) \right) dt = \mathrm{E}  \int_0^T \left(
v(X^{x, \varepsilon}_t)  -  \overline{v}(z^{x,\varepsilon}_{t})  \right) f'(z^{x,\varepsilon}_{t}) dt.
\]
It easily follows from the explicit construction of $X^{x, \varepsilon}_t$ (Section~\ref{fssy}) that the expression in the right hand
side tends to zero uniformly in $x$. 
\qed
\\

 Now let us consider the original situation with two ergodic classes for the Markov chain when $z \geq 0$. Recall that $S$ is now a graph with three semi-infinite
edges, $I_0, I_1$, and $I_2$, with the common vertex $O$. The process $Y^y_t$ on $S$ has been defined in Section~\ref{fssy} (the case $\varkappa = 0$). 
The motion is deterministic on each of the edges, while the behavior at $O$ is random - the process proceeds to $I_1$ or $I_2$ with the prescribed probabilities 
$p_1$ and $p_2$,
respectively.

Let $h$ be the mapping of $M = \{1,...,n\} \times \mathbb{R}$ to $S$ defined as follows:
\begin{equation} \label{projh}
h(i,z) = \left\{ \begin{array}{ll}
                 (0,z), ~~~~~~ z \leq 0,\\
                 (1,z), ~~~~~~ i \in R_1,~z \geq 0,\\
								 (2,z), ~~~~~~ i \in R_2,~z \geq 0.
            \end{array}
            \right.
\end{equation}
\begin{theorem} Suppose that $\varkappa = 0$ and that the assumptions made in Section~\ref{fssy} are satisfied (in particular, the Markov chain $\Xi^z_t$ has two ergodic classes for each $z \geq 0$).
For each $x \in M$, the measures induced by the processes $Y^{x, \varepsilon}_t = h(X^{x, \varepsilon}_t)$ on $S$ converge weakly, as $\varepsilon \downarrow 0$, to
the measure induced by the process $Y^{h(x)}_t$. 
\end{theorem}
\proof
Lemma~\ref{cnvl} is not directly applicable now because the semigroup that corresponds to the process $Y^y_t$ does not preserve $C_0(S)$. 
However, outside of an arbitrarily small neighborhood of the set $h^{-1}(O)$, the limiting motion of $Y^{x, \varepsilon}_t$ is given by $Y^{h(x)}_t$, as follows from Theorem~\ref{saba}. To complete the proof, we need to show that if $X^{x, \varepsilon}_t$ starts slightly to the left of $h^{-1}(O)$, then it quickly moves to the right  of $h^{-1}(O)$ and $\xi^{x, \varepsilon}_t$ ends up in the first ergodic class with probability close to $p_1$.

 More precisely, let
\[
\tau^{x,\varepsilon}_\delta = \inf\{t \geq 0: z^{x,\varepsilon}_t = \delta\}.
\]
It is sufficient to show that for each $\eta > 0$ there is $\delta_0 > 0$ such that each $\delta \in (0,\delta_0]$ there is $\varepsilon_0 > 0$ such that
for $\varepsilon \in (0, \varepsilon_0]$, we have
\begin{equation} \label{mineqa}
\mathrm{E} \tau^{x,\varepsilon}_\delta < \eta,
\end{equation}
\begin{equation} \label{mineq}
|\mathrm{P} (\xi^{x,\varepsilon}_{\tau^{x,\varepsilon}_\delta} \in R_1) - p_1| < \eta,
\end{equation}
whenever $x = (i, -\delta)$. From the explicit construction of $X^{x, \varepsilon}_t$ (Section~\ref{fssy}), it is clear
that 
$z^{x,\varepsilon}_t$ increases, while on $[-\delta, \delta]$, with the speed that is bounded from below by $\inf_{i, z \in [-\delta, \delta]} v(i,z) > 0$. 
This implies (\ref{mineqa}). 
To prove (\ref{mineq}), we define $f_\varepsilon(i,z)$, $z \in [-\delta, \delta]$, as the solution of the system of ODEs 
\[
\frac{d f_\varepsilon(i, z)}{d z}  =  \frac{(v(i,z))^{-1}}{\varepsilon} \left(Q_i( z) f_\varepsilon(i,z) - \sum_{j \neq i} q_{ij}(z) f_\varepsilon(j,z) \right)
\]
with the terminal condition 
\[
f_\varepsilon(i,\delta) = \overline{e}_i := \left\{ \begin{array}{ll}
                 1, ~~~~~~ i \in R_1,\\
                 0, ~~~~~~ i \in R_2.
            \end{array}
            \right.
\]					
We extend $f_\varepsilon$ to be defined on $M$ so that  $f_\varepsilon \in \mathcal{D}(A^\varepsilon)$. Observe that, by construction, $A^\varepsilon f_\varepsilon
(i,z) = 0$ when $z \in [-\delta, \delta]$. Therefore, applying (\ref{sttime}) with $\tau  = \tau^{x,\varepsilon}_\delta$ and $x = (i, -\delta)$, we obtain
\[
\mathrm{P} (\xi^{x,\varepsilon}_{\tau^{x,\varepsilon}_\delta} \in R_1) = f_\varepsilon(i, -\delta).
\]
Thus it remains to analyze the asymptotics of the solution to the ODE. Let $N(z)$ be the matrix, whose diagonal elements
are $N_{ii}(z) = -(v(i,z))^{-1} Q_i( z) $ and off-diagonal elements are $N_{ij}(z) =  (v(i,z))^{-1} q_{ij}(z)$. Let 
\[
N^\delta = \frac{1}{2 \delta} \int_{-\delta}^\delta N(z) dz. 
\]
Solving the linear ODE, we get
\[
f_\varepsilon(\cdot, -\delta) = \exp(\frac{2\delta}{\varepsilon} N^\delta) \overline{e}.
\]
 When $\delta$ is small, $N^\delta$ is a small
perturbation of the matrix $N(0)$. Namely, let
\[
H^\delta = N^\delta - N(0).
\]
All the entries of $H^\delta$ tend to zero when $\delta \downarrow 0$. Observe that all the off-diagonal entries of $N^\delta$ are positive for each $\delta$,
and the sum of elements in each row is equal to zero. Therefore, zero is the simple eigenvalue of $N^\delta$ with the
right eigenvector equal to $e = (1,...,1)^T$, the real parts of the other eigenvalues are negative.

Let $\Pi^\delta_e(\overline{e})$ be the projection of $\overline{e}$ onto $e$ along the space spanned by the remaining eigenvectors (and generalized eigenvectors) 
of the matrix $N^\delta$. Then
\[
\lim_{\varepsilon \downarrow 0} f_\varepsilon(i, -\delta)  = (\Pi^\delta_e(\overline{e}))_i 
\]
for each $i$, and it remains to show that $(\Pi^\delta_e(\overline{e}))_i $ (which does not depend on $i$) is close to $p_1$ for small $\delta$. 

Observe that zero is the top eigenvalue of $N(0)$ with two linearly independent right eigenvectors $e$ and $\overline{e}$ and two linearly independent left
eigenvectors: 
\[
\pi^1_i =  \left\{ \begin{array}{ll}
\pi_i v_i, ~~~~~ i \in R_1,\\
                 0, ~~~~~~~~ i \in R_2,
            \end{array}
            \right.
\]	
\[
\pi^2_i =  \left\{ \begin{array}{ll}
0, ~~~~~~~~ i \in R_1,\\
                 \pi_i v_i, ~~~~~ i \in R_2,
            \end{array}
            \right.
\]		
where $v_i = v(i,0)$. 
Let $\lambda_1^\delta < 0$ be the eigenvalue of $N^\delta$ with the second-largest real part (the top eigenvalue is zero). It is determined uniquely for
small $\delta$. Let $g^\delta$ be the corresponding right eigenvector (determined up to a constant factor). 
\begin{lemma} The vector $g^\delta$ can be represented as
\begin{equation} \label{decompo}
g^\delta = e + \alpha^\delta \overline{e} + \overline{g}_\delta,
\end{equation}
where $\overline{g}_\delta$ belongs to the space spanned by the eigenvectors (and generalized eigenvectors) of $N(0)$, other than $e$ and $\overline{e}$. The coefficient
$\alpha^\delta$ is
bounded away from zero, and $\overline{g}_\delta$ tends to zero when $\delta \downarrow 0$. 
\end{lemma}
\proof
Let $\bar{i}(\delta)$ be such that $|g^\delta_{\bar{i}(\delta)}| = \max_{1 \leq i \leq n} |g^\delta_i|$. Assume, for now, that $\bar{i}(\delta) 
\in R_1$ for all sufficiently small $\delta$. Then, since $N^\delta$ is a small
perturbation of $N(0)$ and $\lambda_1^\delta \rightarrow 0$ as $\delta \downarrow 0$, 
$N^\delta g^\delta = \lambda_1^\delta g^\delta$ easily implies
that $g^\delta_{i}/g^\delta_{\bar{i}(\delta)} \rightarrow 1$ as $\delta \downarrow 0$ for all $i \in R_1$. 

Let $\tilde{\pi}^\delta$ be the normalized left eigenvector for $N^\delta$ with eigenvalue zero.  From $\tilde{\pi}^\delta N^\delta = 0$ and $N^\delta g^\delta = \lambda_1^\delta g^\delta$ it follows that $\langle g^\delta, \tilde{\pi}^\delta \rangle
= 0$. Let $\tilde{i}(\delta)$ be such that $g^\delta_{\tilde{i}(\delta)} = \max_{i \in R_2} |g^\delta_i|$.
Observe that $\tilde{\pi}^\delta_i \rightarrow \pi^1_i$ for
$i \in R_1$, and $\tilde{\pi}^\delta_i \rightarrow \pi^2_i$ for
$i \in R_2$. Therefore,  
\begin{equation} \label{comab}
c_1 |g^\delta_{\bar{i}(\delta)}| \leq |g^\delta_{\tilde{i}}| \leq c_2 |g^\delta_{\bar{i}(\delta)}|
\end{equation}
 for some
positive constants $c_1$ and $c_2$.  As above,  $g^\delta_{i}/g^\delta_{\tilde{i}(\delta)} \rightarrow 1$ as $\delta \downarrow 0$ for all $i \in R_2$. 
From the facts that $\langle g^\delta, \tilde{\pi}^\delta \rangle = 0$,  $\tilde{\pi}^\delta_i \rightarrow \pi^1_i$ for
$i \in R_1$, and $\tilde{\pi}^\delta_i \rightarrow \pi^2_i$ for
$i \in R_2$, it follows that $g^\delta_i$, $i \in R_1$, are of the opposite sign from $g^\delta_i$, $i \in R_2$.

The vector $g^\delta$ can be represented as a sum of three components, $g^\delta = a^\delta + b^\delta + c^\delta$, where $a^\delta$ is a multiple of $e$, $b^\delta$ is a multiple of $\overline{e}$, and $c^\delta$
is in the space spanned by the eigenvectors (and generalized eigenvectors) of $N(0)$, other than $e$ and $\overline{e}$. Observe that $\|c^\delta\|/\|g^\delta\|
\rightarrow 0$ as $\delta \downarrow 0$ since $e$ and $\overline{e}$ span the eigenspace  corresponding to the top eigenvalue of $N(0)$ and $g^\delta$ belongs
to a small perturbation of that space. Moreover, from
(\ref{comab}) and the fact that $g^\delta_i$, $i \in R_1$, and $g^\delta_i$, $i \in R_2$, are of the opposite sign, it follows that $\|a^\delta\|/\|b^\delta\|$ is bounded from above and below. Therefore, (\ref{decompo}) is possible with $\alpha^\delta$ bounded away from zero and infinity. 

Finally, it remains to note that our condition $\bar{i}(\delta)  \in R_1$ does not lead to any loss of generality.
\qed  
\\
\\
Since $g^\delta$ is the eigenvector of $N^\delta$,  we get
\[
(N(0) + H^\delta)(e + \alpha^\delta \overline{e} + \overline{g}_\delta) = \lambda_1^\delta ( e + \alpha^\delta \overline{e} + \overline{g}_\delta).
\]
Taking the scalar product with $\pi^1$ and $\pi^2$ on both sides and noting that $H^\delta e = 0$, we obtain
\[
\alpha^\delta  \langle H^\delta  \overline{e}, \pi^1 \rangle + \langle H^\delta \overline{g}_\delta , \pi^1 \rangle 
= \lambda_1^\delta \langle e + \alpha^\delta \overline{e} , \pi^1 \rangle, 
\]
\[
\alpha^\delta  \langle H^\delta  \overline{e}, \pi^2 \rangle + \langle H^\delta \overline{g}_\delta , \pi^2 \rangle 
= \lambda_1^\delta \langle e + \alpha^\delta \overline{e} , \pi^2 \rangle. 
\]
Therefore,
\[
\left(\alpha^\delta  \langle H^\delta  \overline{e}, \pi^1 \rangle + \langle H^\delta \overline{g}_\delta , \pi^1 \rangle\right) 
\langle e + \alpha^\delta \overline{e} , \pi^2 \rangle = 
\left( \alpha^\delta  \langle H^\delta  \overline{e}, \pi^2 \rangle + \langle H^\delta \overline{g}_\delta , \pi^2 \rangle  \right)
\langle e + \alpha^\delta \overline{e} , \pi^1 \rangle.
\]
Observe that 
\[
\langle H^\delta \overline{g}_\delta , \pi^1 \rangle = o( \alpha^\delta  \langle H^\delta  \overline{e}, \pi^1 \rangle),~~~ 
\langle H^\delta \overline{g}_\delta , \pi^2 \rangle = o( \alpha^\delta  \langle H^\delta  \overline{e}, \pi^2 \rangle),~~~{\rm as}~~\delta \downarrow 0,
\] 
and therefore, 
\[
\langle H^\delta  \overline{e}, \pi^1 \rangle \langle e + \alpha^\delta \overline{e} , \pi^2 \rangle \sim
\langle H^\delta  \overline{e}, \pi^2 \rangle \langle e + \alpha^\delta \overline{e} , \pi^1 \rangle~~~{\rm as}~~\delta \downarrow 0.
\]
Solving for $\alpha^\delta$ gives
\[
\lim_{\delta \downarrow 0} \alpha^\delta = - 1 - \frac{\left(\sum_{i \in R_1} \sum_{j \in R_2} \overline{q}_{ij} \pi_i \right)
\left( \sum_{i \in R_2} \pi^2_i \right)}{
\left(\sum_{i \in R_2} \sum_{j \in R_1} \overline{q}_{ij} \pi_i \right)
\left( \sum_{i \in R_1} \pi^1_i \right)} = 
- 1 - 
\frac{ \sum_{i \in R_2} \pi^2_i }{\sum_{i \in R_1} \pi^1_i }. 
\]
From (\ref{decompo}), it follows that 
\[
\lim_{\delta \downarrow 0} (\Pi^\delta_e(\overline{e}))_i = - 1/ \lim_{\delta \downarrow 0} \alpha^\delta =
 \frac{ \sum_{i \in R_1} \pi^1_i  }{
 \sum_{i \in R_1} \pi^1_i  +  \sum_{i \in R_2} \pi^2_i  }~,
\]
as required.
\qed

\subsection{The case with diffusion}
Now we consider the fast-slow system $X^{x,\varepsilon}_t = (\xi^{x,\varepsilon}_{t}, z^{x,\varepsilon}_{t})$
defined in Section~\ref{fssy}, with $\varkappa  = 1$. The filtration generated by the process will be denoted by $\mathcal{F}^{x,\varepsilon}_t$.
The process $Y^y_t$ on the graph $S$ is now a diffusion (defined in Section~\ref{limpro} via its generator). The mapping $h$ is the same as in (\ref{projh}). 
\begin{theorem} Suppose that $\varkappa = 1$ and that the assumptions made in Section~\ref{fssy} are satisfied (in particular, the Markov chain $\Xi^z_t$ has two ergodic classes for each $z \geq 0$).
For each $x \in M$, the measures induced by the processes $Y^{x, \varepsilon}_t = h(X^{x, \varepsilon}_t)$ on $S$ converge weakly, as $\varepsilon \downarrow 0$, to
the measure induced by the process $Y^{h(x)}_t$. 
\end{theorem}
\proof
Let $T > 0$, $f \in \mathcal{D}(A)$, and let $K$ be a compact subset of $S$.
It is clear that the family of measures on $C([0, \infty), S)$ induced by the processes $Y^{x, \varepsilon}_t$, $\varepsilon > 0$, is tight for each $x \in M$.
Thus, by Lemma~\ref{cnvl}, it is sufficient to prove that, given $\eta > 0$, we have
\[
	|\mathrm{E} (f( Y^{x, \varepsilon}_T ) - f (  Y^{x, \varepsilon}_0) - \int_0^T A f (  Y^{x, \varepsilon}_t) dt) |  \leq \eta,
	\]
for all $x \in h^{-1} (K)$ and all sufficiently small $\varepsilon$. 

Let us define two sequences of stopping times: 
\[
\sigma^{x,\varepsilon}_0 = 0;~~~\tau^{x,\varepsilon}_n = \inf\{t \geq \sigma_{n-1}: z^{x,\varepsilon}_{t} = 0\},~n \geq 1;~~~ 
\sigma^{x,\varepsilon}_n = \inf\{t \geq \tau_{n}: |z^{x,\varepsilon}_{t}| = \delta \},~n \geq 1,
\]
where $\delta > 0$ will be selected later. Then
\[
\mathrm{E} \left(f( Y^{x, \varepsilon}_T ) - f (  Y^{x, \varepsilon}_0) - \int_0^T A f (  Y^{x, \varepsilon}_t) dt\right)=
\]
\begin{equation} \label{sums}
 \mathrm{E}  \sum_{n=1}^\infty \left(f( Y^{x, \varepsilon}_{\tau^{x,\varepsilon}_n \wedge T} ) - f (  Y^{x, \varepsilon}_{\sigma^{x,\varepsilon}_{n-1} \wedge T}) - 
\int_{\sigma^{x,\varepsilon}_{n-1} \wedge T}^{\tau^{x,\varepsilon}_n \wedge T} A f (  Y^{x, \varepsilon}_t) dt\right) +
\end{equation}
\[
 \mathrm{E}  \sum_{n=1}^\infty \left(f( Y^{x, \varepsilon}_{\sigma^{x,\varepsilon}_n \wedge T} ) - f (  Y^{x, \varepsilon}_{\tau^{x,\varepsilon}_{n} \wedge T}) - 
\int_{\tau^{x,\varepsilon}_{n} \wedge T}^{\sigma^{x,\varepsilon}_n \wedge T} A f (  Y^{x, \varepsilon}_t) dt\right).
\]
In order to control the number of terms in the sums above, we'll need the following lemma.
\begin{lemma} \label{lennk}
There is $c > 0$ such that, for all sufficiently small $\delta$,
\begin{equation} \label{oll1}
\mathrm{P}(\sigma^{x,\varepsilon}_{n} \leq T) \leq \exp(-c \delta n),~~~x \in M,~~n \geq 2.
\end{equation}
\end{lemma}
\proof 
Let $A_t$ be an auxiliary diffusion process, $dA_t = a dt + d W_t$, $A_0 = -\delta$, where
$a = \sup_{i,z}|v(i,z)|$. Let $\tilde{\tau} = \inf\{t: A_t = 0\}$. Then 
 $\mathrm{P}(\tilde{\tau} \leq T) \leq \exp(-c \delta)$ for some $c > 0$. If $\tilde{\tau}_k$, $k \geq 1$, is a sequence of
independent random variables distributed as $\tilde{\tau}$, then
\begin{equation} \label{tyt}
\mathrm{P}(\tilde{\tau}_1 + ... + \tilde{\tau}_n \leq T) \leq \exp(-c \delta n).
\end{equation}
From the definition of the stopping times and the process $X^{x,\varepsilon}_t$ it follows that 
\[
\mathrm{P}(\tau^{x,\varepsilon}_n - \sigma^{x,\varepsilon}_{n-1}
> s | \mathcal{F}^{x,\varepsilon}_{\sigma^{x,\varepsilon}_{n-1}} ) \geq  \mathrm{P}(\tilde{\tau} > s)
\]
for each $n \geq 2$ and $s \geq 0$. Therefore, estimate (\ref{oll1}), with $\tau^{x,\varepsilon}_{n+1}$ instead of $\sigma^{x,\varepsilon}_{n}$,
follows from (\ref{tyt}) and the strong Markov property. Thus, the original formula  (\ref{oll1}) also holds, with a different constant $c$. \qed
\\

Let
\[
\alpha(x,n) = \mathrm{E}  \left(f( Y^{x, \varepsilon}_{\tau^{x,\varepsilon}_n \wedge T} ) - f (  Y^{x, \varepsilon}_{\sigma^{x,\varepsilon}_{n-1} \wedge T}) - 
\int_{\sigma^{x,\varepsilon}_{n-1} \wedge T}^{\tau^{x,\varepsilon}_n \wedge T} A f (  Y^{x, \varepsilon}_t) dt  | 
\mathcal{F}^{x,\varepsilon}_{\sigma^{x,\varepsilon}_{n-1} \wedge T} \right).
\]
Observe that
\[
\lim_{\varepsilon \downarrow 0} \sup_{x \in h^{-1} (K)} \sup_{n \geq 1}  |\alpha(x,n) | = 0
\]
 uniformly in all the realizations of the randomness (which is present since we are taking the conditional expectation). 
 This is a standard averaging result for the fast-slow system in the case of a single invariant measure for the fast motion. 
It easily follows from the explicit construction of $X^{x, \varepsilon}_t$. Therefore, for the first expectation in (\ref{sums}), by Lemma~\ref{lennk}, we get
\[
|\mathrm{E}  \sum_{n=1}^\infty \left(f( Y^{x, \varepsilon}_{\tau^{x,\varepsilon}_n \wedge T} ) - f (  Y^{x, \varepsilon}_{\sigma^{x,\varepsilon}_{n-1} \wedge T}) - 
\int_{\sigma^{x,\varepsilon}_{n-1} \wedge T}^{\tau^{x,\varepsilon}_n \wedge T} A f (  Y^{x, \varepsilon}_t) dt\right)| 
\leq
\]
\[
 \leq \sum_{n=1}^\infty |\alpha(x,n) | \mathrm{P}(\sigma^{x,\varepsilon}_{n-1} \leq T) \rightarrow 0~~~{\rm as}~~\varepsilon \downarrow 0,
\]
uniformly in $x \in h^{-1} (K)$.

Next, observe that 
\[
| \mathrm{E}  \left({\sigma^{x,\varepsilon}_n \wedge T} -  {\tau^{x,\varepsilon}_n \wedge T} | 
\mathcal{F}^{x,\varepsilon}_{\tau^{x,\varepsilon}_n \wedge T} \right) | \leq C \delta^2
\]
for some constant $C$ and all $x \in M$, $n \geq 1$. This follows from the fact that the process $z^{x,\varepsilon}_{t}$ is a Brownian motion with a bounded variable drift, and
the expectation of its exit time from the $\delta$-neighborhood of the origin is estimated from above by $C \delta^2$.  Therefore, 
\[
| \mathrm{E}  \sum_{n=1}^\infty 
\int_{\tau^{x,\varepsilon}_{n} \wedge T}^{\sigma^{x,\varepsilon}_n \wedge T} A f (  Y^{x, \varepsilon}_t) dt| \leq  C \delta^2 \sup|Af| \sum_{n=1}^\infty 
 \mathrm{P}(\tau^{x,\varepsilon}_{n} \leq T). 
\]
By  Lemma~\ref{lennk}, since $\tau^{x,\varepsilon}_{n} \geq \sigma^{x,\varepsilon}_{n-1}$, the right hand side does not exceed $K \delta$ for some constant $K$. 
This is smaller than $\eta/2$ for all sufficiently small $\delta$. Thus it remains to show that there is $\delta > 0$ such that
\[
 | \mathrm{E}  \sum_{n=1}^\infty \left(f( Y^{x, \varepsilon}_{\sigma^{x,\varepsilon}_n \wedge T} ) - 
f (  Y^{x, \varepsilon}_{\tau^{x,\varepsilon}_{n} \wedge T}) \right)|  < \eta/2
\]
for all sufficiently small $\varepsilon$. Observe that
\[
| \mathrm{E}  \sum_{n=1}^\infty \left(f( Y^{x, \varepsilon}_{\sigma^{x,\varepsilon}_n} ) - 
f (  Y^{x, \varepsilon}_{\sigma^{x,\varepsilon}_n \wedge T}) \right) \chi_\{ \tau^{x,\varepsilon}_{n} \leq T \}|  \leq \sup|f(l_1,z_1) - f(l_2,z_2)| < \eta/4
\]
for all sufficiently small $\delta$,
where the supremum is taken over all $l_1, l_2$ and $z_1, z_2$ such that $|z_1|, |z_2| \leq \delta$. Therefore,
\[
| \mathrm{E}  \sum_{n=1}^\infty \left(f( Y^{x, \varepsilon}_{\sigma^{x,\varepsilon}_n \wedge T} ) - 
f (  Y^{x, \varepsilon}_{\tau^{x,\varepsilon}_{n} \wedge T}) \right)| \leq \eta/4 +  \sup_{x: h(x) = O} \mathrm{E}   
\left(f( Y^{x, \varepsilon}_{\sigma} ) - 
f (  Y^{x, \varepsilon}_0) \right) \sum_{n=1}^\infty 
 \mathrm{P}(\tau^{x,\varepsilon}_{n} \leq T),
\]
where $\sigma = \sigma^{x,\varepsilon} = \inf\{t \geq 0: |z^{x,\varepsilon}_{t}| = \delta \}$. By  Lemma~\ref{lennk}, since 
$\tau^{x,\varepsilon}_{n} \geq \sigma^{x,\varepsilon}_{n-1}$, the sum in the right hand side can be estimated from above by $K/\delta$ for some $K$, and it remains
to show that $\sup_{x: h(x)  = O} \mathrm{E}   \left(f( Y^{x, \varepsilon}_{\sigma} ) - 
f (  Y^{x, \varepsilon}_0) \right) /\delta$ can be made arbitrarily small for some $\delta$ and all sufficiently small $\varepsilon$. Since $f(l,z)$ is differentiable
in $z$ at $z = 0$ along each edge (one-sided derivatives exist), and the relation between the derivatives is given by (\ref{relder}), the result follows from
the following lemma.
\begin{lemma} \label{cllee}
For each $\eta > 0$, for all sufficiently small $\delta > 0$, 
\[
|\mathrm{P}(\xi^{x,\varepsilon}_{\sigma} \in R_1, z^{x,\varepsilon}_{\sigma} = \delta   ) - \frac{\overline{\pi}_1}{2} | \leq \eta,~~~~
|\mathrm{P}(\xi^{x,\varepsilon}_{\sigma} \in R_2, z^{x,\varepsilon}_{\sigma} = \delta   ) - \frac{\overline{\pi}_2}{2} | \leq \eta 
\]
for each $x$ such that $h(x) = O$ and  all sufficiently small $\varepsilon$ (depending on $\delta$).
\end{lemma}
\proof
Consider an auxiliary process  $\tilde{X}^{x,\varepsilon}_t = (\tilde{\xi}^{x,\varepsilon}_{t}, \tilde{z}^{x,\varepsilon}_{t})$ that is defined
the same way as  $X^{x,\varepsilon}_t$, but with $v(i, \cdot) \equiv 0$ for each $i$. The corresponding stopping time will
be denoted by $\tilde{\sigma}$. Let $\tilde{\mu}_{t_0}$ and $\mu_{t_0}$ be the measures on the space of 
RCLL functions from $[0,t_0]$ to $M$ induced by the processes $\tilde{X}^{x,\varepsilon}_t$ and ${X}^{x,\varepsilon}_t$, respectively.
By the Girsanov theorem,  $\tilde{\mu}_{t_0}$ and $\mu_{t_0}$ are mutually absolutely continuous. Moreover, for each $\eta > 0$, 
for all sufficiently small $t_0$ and $\varepsilon$, we have $\tilde{\mu}_{t_0} (1-\eta \leq p_{t_0} \leq 1 + \eta) \geq 1 - \eta$ for each $x = (i, 0)$, where
$p_{t_0}$ is the density of $\mu_{t_0}$ with respect to $\tilde{\mu}_{t_0}$. Since $\mathrm{P}(\tilde{\sigma} \leq t_0) \rightarrow 1$ as $\delta \downarrow 0$,  
 we have, for all sufficiently small $\delta$ and $\varepsilon$ and $x = (i, 0)$,
\[
\lim_{\delta \downarrow 0}  |\mathrm{P}(\xi^{x,\varepsilon}_{\sigma} \in R_1, z^{x,\varepsilon}_{\sigma} = \delta   )  - 
\mathrm{P}(\tilde{\xi}^{x,\varepsilon}_{\tilde{\sigma}} \in R_1, \tilde{z}^{x,\varepsilon}_{\tilde{\sigma}} = \delta   )| \leq \eta/2.
\]
Similarly,
\[
\lim_{\delta \downarrow 0}  |\mathrm{P}(\xi^{x,\varepsilon}_{\sigma} \in R_2, z^{x,\varepsilon}_{\sigma} = \delta   )  - 
\mathrm{P}(\tilde{\xi}^{x,\varepsilon}_{\tilde{\sigma}} \in R_2, \tilde{z}^{x,\varepsilon}_{\tilde{\sigma}} = \delta   )| \leq \eta/2.
\]
Thus it is sufficient to prove Lemma~\ref{cllee} in the case when there is no drift term. 

Next, we need the following observation about time-inhomogeneous Markov processes. 
Recall that the time-homogeneous Markov chain with transition rates $q_{ij}(z)$, $z < 0$,
has a unique invariant distribution $\mu_i(z)$, $1 \leq i \leq n$, $z \in  \mathbb{R}$. Moreover, when $t \rightarrow \infty$, the distribution of $\Xi^z_t$
is close to the invariant distribution, which, in turn is close to $\pi_i$, $1 \leq i \leq n$, if $|z|$ is small. A similar statement can be made 
about time-inhomogeneous processes. Namely, let $\eta > 0$. It is not difficult to show that there is $\delta > 0$ with the following property: if $|\tilde{z}(t)| \leq \delta$ for $t \leq t_0$ and 
if there is $\delta_0 > 0$ such that $\lambda(t: \tilde{z}(t) \leq - \delta_0) \rightarrow \infty$ as $t_0 \rightarrow \infty$, then
\[
\lim_{t_0 \rightarrow \infty} | \mathrm{P}(\tilde{\Xi}^{\tilde{z}}_{t_0} - i) - \pi_i| \leq \eta, 
\]
where $\lambda$ is the Lebesgue measure on the real line and $\tilde{\Xi}^{\tilde{z}}_t$ is a time-inhomogeneous Markov process with transition rates at time
$t$ given by $q_{ij}(\tilde{z}(t))$.

To complete the proof of Lemma~\ref{cllee} in the case when there is no drift term, we condition the evolution of the fast component on the realization of
the Brownian motion and obtain that the above argument is applicable for almost every realization of the Brownian motion (after rescaling the time by $1/\varepsilon$). 
\qed
\\

As we discussed above, Lemma~\ref{cllee} completes the proof of the theorem. 
\qed
%
\\\\
\noindent {\bf \large Acknowledgments}:  While working on this
article, L. Koralov was supported by the ARO grant W911NF1710419 and by the University of Maryland Research and Scholarship Award. 
\\
\\


\begin{thebibliography}{999999}


%










\bibitem{Dyn} Dynkin E. B., {\it Markov Processes}, Springer-Velag, Berlin, Heidelberg, New York, 1965.


\bibitem{EK86} Ethier S. N.,  Kurtz T. G, {\it  Markov processes: characterization and convergence},  Wiley Series
in Probability and Mathematical Statistics: Probability and Mathematical Statistics. John Wiley and Sons,
Inc., New York, 1986.

\bibitem{F85} Freidlin M.I., {\it Functional Integration and
Partial Differential Equations}, Princeton University Press, 1985.

\bibitem{Fp1} Freidlin M.I., {\it Thermostat-like perturbations of an oscillator}, J. Stat. Phys. 164 (2016), no. 1, pp. 130--141.

\bibitem{Fp2} Freidlin M.I., {\it On stochastic perturbations of dynamical systems with a ``rough" symmetry. Hierarchy of Markov chains}, 
J. Stat. Phys. 157 (2014), no. 6, 1031--1045. 

\bibitem{FW} Freidlin M. I., Wentzell A. D., {\it Random
Perturbations of Dynamical Systems}, Springer 2012.

\bibitem{KPST} Korolyuk V. S.,  Portenko N. I,  Skorokhod A. V.,  Turbin A. F., {\it Handbook on probability theory and mathematical statistics}, Nauka, 1985,
(in Russian). 









%
\bibitem{Mandl} Mandl P., {\it Analytical Treatment of One-dimensional Markov Processes}, Springer-Verlag, 1968.






\end{thebibliography}
\end{document}